\documentclass[a4paper, oneside]{article}
\usepackage[scale=0.75, centering]{geometry}
\usepackage{pdfsync}
\usepackage{amsmath,amscd,amssymb,amsthm}
\usepackage{graphicx}
\usepackage[english]{babel}   
\usepackage[utf8]{inputenc}
\usepackage[T1]{fontenc}
\usepackage{lmodern}
\usepackage{array}
\usepackage{xspace}
\usepackage{hyperref}
\newcommand{\eq}[1]{\begin{equation} #1 \end{equation}}
\newcommand{\dofs}{DoFs\@\xspace}
\newcommand{\Real}{{\mathbb{R}}}
\newcommand{\vect}[1]{\underline{#1}}
\newcommand{\norma}[1]{\vect{\nu}_{#1}}
\newcommand{\pOm}{\partial\Omega}
\newcommand{\Vit}{u}
\newcommand{\Vex}{\vect{\Vit}}
\newcommand{\Vort}{\omega}
\newcommand{\Wex}{\vect{\Vort}}
\newcommand{\Pex}{p}
\newcommand{\Foex}{\vect{f}}
\newcommand{\vtxh}{\mathrm{V}}
\newcommand{\edgeh}{\mathrm{E}}
\newcommand{\faceh}{\mathrm{F}}
\newcommand{\cellh}{\mathrm{C}}
\DeclareMathOperator{\grd}{\vect{grad}}
\DeclareMathOperator{\divm}{div}
\DeclareMathOperator{\curl}{\vect{curl}}
\newcommand{\HypO}{($\mathbf{H_{\Omega}}$)\@\xspace}

\begin{document}

\title{Analysis of Compatible Discrete Operator Schemes for the Stokes Equations on Polyhedral Meshes}

\author{
\begin{tabular}[t]{c@{\extracolsep{4em}}c}
Jerome Bonelle & Alexandre Ern \\
\begin{small}
EDF R\&D
\end{small}
&
\begin{small}
Universit\'e Paris-Est, CERMICS
\end{small} \\
\begin{small}
 6, quai Watier, BP 49
\end{small}
&
\begin{small}
Ecole des Ponts ParisTech
\end{small} \\
\begin{small}
78401 Chatou cedex
\end{small}
&
\begin{small}
77455 Marne la Vall\'ee Cedex 2, France
\end{small} \\
\begin{small}
jerome.bonelle@edf.fr
\end{small}
&
\begin{small}
ern@cermics.enpc.fr
\end{small}
\end{tabular}}

\maketitle

\begin{abstract}
{
Compatible Discrete Operator schemes preserve basic properties of the continuous model at the discrete level. They combine discrete differential operators that discretize exactly topological laws and discrete Hodge operators that approximate constitutive relations. We devise and analyze two families of such schemes for the Stokes equations in curl formulation, with the pressure degrees of freedom located at either mesh vertices or cells. The schemes ensure local mass and momentum conservation. We prove discrete stability by establishing novel discrete Poincar\'e inequalities. Using commutators related to the consistency error, we derive error estimates with first-order convergence rates for smooth solutions. 
We analyze two strategies for discretizing the external load, so as to deliver tight error estimates when the external load has a large irrotational or divergence-free part.
Finally, numerical results are presented on three-dimensional polyhedral meshes.
The detailed material is available from
\url{http://hal.archives-ouvertes.fr/hal-00939164}

}
\end{abstract}

\section{Introduction}
\label{sec:intro}

Compatible Discrete Operator (CDO) schemes belong to the broad class of compatible or mimetic schemes, which preserve basic properties of the continuous model at the discrete level; see \cite{AnBHK:12, ArFaW:10, BocHy:05, BrBuL:09, BrLiS:05, CleWe:01, DeHLM:05, Gerri:12, Matti:00, Perot:11, Teixe:13, Tonti:01} and references therein.
Following the seminal ideas of~\cite{Tonti:75} and~\cite{Bossa:00}, the degrees of freedom (\dofs) are defined using de Rham maps, and their localization results from the physical nature of the fields.
Moreover, a distinction is operated between topological laws (that are discretized exactly) and constitutive relations (that are approximated). 
CDO schemes are formulated using discrete differential operators for the topological laws and discrete Hodge operators for
the constitutive relations.
The discrete differential operators produce a cochain complex, they commute with the de Rham maps, and discrete adjunction properties hold between these operators.
The discrete Hodge operator is the key operator in the CDO framework.
The design of this operator is not unique, and each design leads to a specific scheme; see~\cite{TaKeB:99, Hiptm:01, BonEr:14}.
CDO schemes involve two meshes: a primal mesh (which is the only one seen by the end-user) and a dual mesh.
The discrete Hodge operator links \dofs defined on the primal mesh to \dofs defined on the dual mesh.
In~\cite{BonEr:14}, CDO schemes have been analyzed for elliptic problems on polyhedral meshes.

The Stokes equations model flows of incompressible and viscous fluids where the advective inertial forces are negligible with respect to the viscous forces.
In this paper, we focus on the stationary Stokes equations posed on an open, bounded and connected domain $\Omega \subset \Real^3$ with boundary $\pOm$ and outward normal $\norma{\pOm}$.
Our starting point is to formulate the viscous stresses in the momentum balance using the curl operator.
This way, all the terms in the Stokes equations can be interpreted using scalar-valued differential forms.
We analyze two formulations.
The first one, hereafter called 2-field curl formulation, takes the form
\eq{\label{eq:Stokes.curl.twofields}\left\lbrace
\begin{aligned}
	&\curl(\curl(\Vex)) + \grd(\Pex) &={}& \Foex, &\qquad \text{in $\Omega$}, \\
	&\divm(\Vex) &={}& 0, &\qquad \text{in $\Omega$},
\end{aligned}\right.}
where $\Pex$ is the pressure, $\Vex$ the velocity and $\Foex$ the external load.
Introducing the vorticity  $\Wex := \curl\Vex$, the second formulation, hereafter called 3-field curl formulation (also called Velocity-Vorticity-Pressure formulation in the literature), takes the form
\eq{\label{eq:Stokes.curl.threefields}\left\lbrace
\begin{aligned}
	&-\Wex  + \curl(\Vex) &={}& \vect{0}, &\qquad \text{in $\Omega$}, \\
	&\curl(\Wex) + \grd(\Pex) &={}& \Foex, &\qquad \text{in $\Omega$}, \\
	&\divm(\Vex) &={}& 0, &\qquad \text{in $\Omega$}.
\end{aligned}\right.}
Essential and natural boundary conditions (BCs) can be considered for both formulations.
The first set of BCs enforces the value of the normal component of the velocity $\Vex\cdot\norma{\pOm}$ and that of the tangential components of the vorticity $\Wex\times\norma{\pOm}$ at the boundary. 
These BCs are natural for~\eqref{eq:Stokes.curl.twofields} and essential for~\eqref{eq:Stokes.curl.threefields}.
As the pressure is then determined up to an additive constant, the additional requirement of $\Pex$ having zero mean-value is typically added. 
The second set of BCs enforces the value of the tangential components of the velocity $\Vex\times\norma{\pOm}$ and the value of the pressure at the boundary. 
These BCs are essential for~\eqref{eq:Stokes.curl.twofields} and natural for~\eqref{eq:Stokes.curl.threefields}.

Mimetic or compatible schemes for the Stokes equations using either the 2-field or the 3-field curl formulation with the above BCs have already been investigated.
Most of the work dedicated to the Stokes equations in curl formulations addresses the 3-field curl formulation~\eqref{eq:Stokes.curl.threefields}.
Based on the seminal work of~\cite{Nedel:82}, \cite{Duboi:92, Duboi:02}
first analyzed this formulation in the context of finite elements (FE), and
\cite{BerCh:06} in the context of spectral discretization.
More recently, \cite{KreGe:13} proposed a scheme based on the Mimetic Spectral Element method in 2D and 3D for general quadrilateral/hexahedral meshes, and \cite{DelOm:13} designed a Discrete Duality Finite Volume (DDFV) scheme for general 2D meshes.
Concerning the 2-field curl formulation, only the papers
of~\cite{BraLe:94} and of~\cite{AbFiS:12} have addressed this
formulation on simplicial meshes; see also~\cite{PerNa:03}
and~\cite{EyFuL:11} for 
staggered schemes on triangular Delaunay meshes. 

In the present work, we devise and analyze CDO schemes for the Stokes problem in both 2-field and 3-field curl formulations.
The CDO schemes involve two Hodge operators, one linking the velocity (seen as a circulation) to the mass flux and the other linking the vorticity to the viscous stress.
One key feature of the present schemes is that they ensure local mass and momentum conservation on polyhedral meshes.
We prove discrete stability by establishing novel discrete Poincar\'e inequalities in the CDO framework.
The present schemes do not need any stabilization.
Moreover, using commutators related to the consistency error as in~\cite{Bossa:00, Hiptm:01, BonEr:14}, we derive \textit{a priori} error estimates and establish first-order error estimates for smooth solutions.
In addition, we show how the present CDO framework can deal with the practically important issue of discretizing the external load, so as to deliver tight error estimates when the external load is expected to have a large irrotational or divergence-free part (see \cite{Linke:14} for a related work).

The more classical formulation of the Stokes equations uses the vector
Laplacian of the velocity in the momentum balance. Various schemes were
proposed to discretize this formulation on polygonal or polyhedral meshes,
including Mixed Finite Volumes by~\cite{DroEy:06}, Mimetic Finite
Differences (MFD) by~\cite{BeGLM:09, BeLiM:10}, DDFV by~\cite{KreMa:12},
and an extension of the Crouzeix--Raviart finite element by~\cite{DPiLe:14}.
\cite{FalNe:13} investigated a scheme for triangular meshes within the Finite Element Exterior Calculus (FEEC) framework.
Discretizing the vector Laplacian with CDO schemes is the subject of ongoing work.

This paper is organized as follows.
In Section~2, we briefly recall the CDO framework.
Then, we investigate CDO schemes for the Stokes equations in the curl formulations~\eqref{eq:Stokes.curl.twofields} and~\eqref{eq:Stokes.curl.threefields}.
Since the pressure is seen as a potential, its \dofs are located at primal or dual mesh vertices.
The former case, treated in Section~3, hinges on the 2-field curl formulation~\eqref{eq:Stokes.curl.twofields} leading to vertex-based pressure schemes.
The latter case, treated in Section~4, hinges on the 3-field curl formulation~\eqref{eq:Stokes.curl.threefields} leading to cell-based pressure schemes (since primal cells are in one-to-one correspondence with dual mesh vertices). The vertex-based pressure schemes are, to our knowledge, the first of this class on polyhedral meshes. The cell-based pressure schemes share common features with the recent Mimetic Spectral Element schemes of~\cite{KreGe:13}; the present schemes can be deployed on polyhedral meshes and offer two strategies for discretizing the external load.
Finally, we present numerical results in Section~5.
For simplicity, we often assume in what follows that $\Omega$ is simply
connected and that its boundary $\pOm$ is connected. Whenever needed, we refer to this assumption as \HypO.

\vspace*{2cm}
\begin{center}
\begin{tabular}{c}
\hline
	\\
	The detailed material is available from \url{http://hal.archives-ouvertes.fr/hal-00939164} \\
	\\
\hline
\end{tabular}
\end{center}
\vspace*{1cm}

\section{Conclusions}
\label{sec:conclusion}

In this work, we have analyzed CDO schemes for the Stokes equations on three-dimensional polyhedral meshes.
The distinction between primal and dual meshes enabled us to devise vertex-based and cell-based pressure schemes.
Vertex-based pressure schemes lead to an algebraic system of size $(\#\vtxh + \#\edgeh)$ with two unknowns, the pressure located at primal vertices and the velocity at primal edges. Cell-based pressure schemes lead to a system of size $(\#\edgeh + \#\faceh + \#\cellh)$ with three unknowns, the pressure located at primal cells, the mass flux at primal faces, and the viscous stress circulation at primal edges. For both schemes, two discrete Hodge operators, related to the mass density and the viscosity, are used; as for elliptic problems, these operators must satisfy a stability and a consistency property. Both schemes conserve mass and momentum, the vertex-based ones at dual cells and dual faces, respectively, and the cell-based ones at primal cells and dual edges, respectively. Finally, both schemes can be deployed with two possible load discretizations, so as to handle a large irrotational or divergence-free part of the load. Finally, various tracks are worth pursuing in future work, including the use of hybridization techniques for saddle-point problems and the treatment of essential (and more general) BCs as well as complex topologies for the domain $\Omega$.

%
%

\end{document}